\title{Polynomial Representation for the Expected Length of Minimal Spanning Trees\footnote{This research was supported by NSF grant \#0649068 funding the WiVaM REU-RET in Mathematics}}

\author{Jared Nishikawa \thanks{Willamette University}
        \and Peter T. Otto $^\dag$
        \and Colin Starr $^\dag$
        }
\documentclass{article}
\usepackage{graphicx}
\usepackage{epsfig}
\usepackage{amsmath, amsthm, amssymb}

\newtheorem{defn}{Definition}[section]

\newtheorem{example}[defn]{Example}

\newtheorem{lemma}[defn]{Lemma}

\newtheorem{proposition}[defn]{Proposition}

\newtheorem{remark}[defn]{Remark}

\newtheorem{thm}[defn]{Theorem}

\newcommand{\noi}{\noindent}

\newcommand{\hsp}{\hspace{.2in}}
\newcommand{\goto}{\rightarrow}

\newcommand{\bea}{\begin{eqnarray}}
\newcommand{\eea}{\end{eqnarray}}
\newcommand{\beas}{\begin{eqnarray*}}
\newcommand{\eeas}{\end{eqnarray*}}
\newcommand{\be}{\begin{equation}}
\newcommand{\ee}{\end{equation}}

\newcommand{\per}{\hspace{-.072in}{\bf .  }}

\newcommand{\brmk}{\begin{remark}\per\begin{em}}
\newcommand{\ermk}{\end{em}\end{remark}}

\newcommand{\bexa}{\begin{example}\per\begin{em}}
\newcommand{\eexa}{\end{em}\end{example}}
\newcommand{\dstyle}{\displaystyle}
\newcommand{\bp}{\begin{proof}}
\newcommand{\ep}{\end{proof}}

\newcommand{\ba}{\begin{array}}
\newcommand{\ea}{\end{array}}

\begin{document}
\newpage
\maketitle
\begin{abstract}
In this paper, we investigate the polynomial integrand of an integral formula that yields the expected length of the minimal spanning tree of a graph whose edges are uniformly distributed over the interval $[0,1]$.  In particular, we derive a general formula for the coefficients of the polynomial and apply it to express the first few coefficients in terms of the structure of the underlying graph; e.g. number of vertices, edges and cycles.

\end{abstract}



\section{Introduction}

In 2002, J.M. Steele \cite{Ste02} derived an integral formula for the expected length of a minimal spanning tree (MST) of a graph with independent edge lengths uniformly distributed over the interval $[0,1]$.  While the formula gives an exact value of the mean length of the MST in terms of the Tutte polynomial of the graph, it yields (at least to us) little intuition of how the MST relates to the structure of the underlying graph.

This provided the motivation for the research project investigated by the Willamette University group of the Willamette Valley REU-RET Consortium for Mathematics Research in the summer of 2008.  The authors of this paper were members of that research group and this paper covers the work that began that summer.

The main result of this paper is a general formula for the coefficients of the polynomial integrand in Steele's formula for the expected length of the MST of a simple, finite, connected graph.  For the first seven coefficients of the polynomial, we prove a surprising result that  expresses the coefficients in terms of features of the underlying graph; e.g. the number of vertices, edges, and cycles.

The remainder of this paper is organized as follows: In Section \ref{formula}, we state Steele's formula, which is written in terms of the Tutte polynomial of the underlying graph.  In Section \ref{integrand}, we investigate the integrand of the formula and prove that it is a polynomial, expressing the coefficients in terms of characteristics of the graph.  We illustrate our results with an example in Section \ref{example} and examine the particular case of the complete graph in Section \ref{complete}.

Throughout this paper, ``graph'' means a finite simple graph.  We adopt the usual notations: $V(G)$ and $E(G)$ are the vertex and edge sets of $G,$ respectively.  The rank of $G$ is denoted by $r(G),$ and $r(G)=|V(G)|-k(G),$ where $k(G)$ denotes the number of connected components of $G.$

\section{Steele's Formula}\label{formula}  

Let $G$ be a graph.  We assign independent random lengths $\xi_e$ with uniform distribution over the interval $[0,1]$ to the edges $e \in E(G)$.  The total length of a minimal spanning tree (MST) of the graph $G$ is denoted by \[ L(G) =\sum_{e \in E(\mbox{{\scriptsize MST}}(G))} \xi_e. \]
We are interested in the expected value of $L(G)$, which we denote by $\mathbb{E}[L(G)]$.

Steele's formula for $\mathbb{E}[L(G)]$ is written in terms of the Tutte polynomial of a graph, which we define next.
\begin{defn}
    \label{defn:tuttepoly}
Let $G$ be a graph, and define $S(G)$ to be the set of spanning subgraphs of $G$; i.e., subgraphs of $G$ with vertex set $V(G)$ and edge set a subset of $E(G)$.  The {\bf Tutte polynomial} of $G$ is defined as follows:
    \[ T(G;x,y)=\sum_{A \in S(G)}^{ }
    (x-1)^{r(G)-r(A)}(y-1)^{|E(A)|-r(A)}. \]
    
\end{defn}
\noi
The Tutte polynomial of a graph encodes much information about the graph, but we will only use the definition above in our analysis and refer the reader to \cite{Bo} for more information.

We will use the following result about the Tutte polynomial in the proof of the main result.  The proof is a straightforward calculation from the definition and so we omit it.

\begin{lemma}
\label{lemma:simpletutte}
Let $G$ be a connected graph with $n$ vertices and $m$ edges.  Then for values of $(x,y)$ satisfying $(x-1)(y-1) = 1$, we have

\medskip

\noi
{\em (a)}  $ {\dstyle T(G; x,y) = (x-1)^{n-1} \left(\frac{x}{x-1}\right)^{m} }$

\medskip

\noi
{\em (b)} ${\dstyle T_x(G; x,y) = (x-1)^{n-2} \left[ \sum_{A \in S(G)} k(A) (y-1)^{|A|} - \left( \frac{x}{x-1} \right)^{m} \right] }$

where $T_x$ denotes the partial derivative of $T$ with respect to $x$.
\end{lemma}

We now state Steele's integral formula for the expected length of the minimal spanning tree that was proved in \cite{Ste02}.

\begin{thm}
\label{thm:main} {\em ({\bf Steele's formula})} \ Let $G$ be a connected graph and $T(G;x,y)$ the Tutte polynomial of $G$.  Then
\be
\label{eqn:steele}
\mathbb{E}[L(G)] = \int_0^1 \frac{1-t}{t} \frac{T_x\left(G; \frac{1}{t},\frac{1}{1-t}\right)}{T\left(G; \frac{1}{t}, \frac{1}{1-t}\right)} \, dt
\ee
\end{thm}

Steele's formula above has been generalized to the case of an arbitrary, but still identical, edge distribution \cite{LZ} and to edge distributions that are not necessarily identically distributed \cite{LZ10}.

\section{Integrand in Steele's Formula}\label{integrand} 

\medskip

\subsection{Polynomial integrand} 

We begin by showing that the integrand in Steele's integral formula is a polynomial of degree less than or equal to the number of edges in the graph.

\begin{thm}
\label{thm:pmt}
Let $G$ be a connected graph with $n$ vertices and $m$ edges.  Then
\[ \mathbb{E}[L(G)] = \int_0^1 p_m(t) \, dt \]
where $p_m(t)$ is a polynomial of degree less than or equal to $m$.
\end{thm}

\bp
For convenience, we let $|A| = |E(A)|$.  By Lemma \ref{lemma:simpletutte}, we have
\beas
\label{eqn:integrand}
\frac{1-t}{t} \frac{T_x\left(G; \frac{1}{t},\frac{1}{1-t}\right)}{T\left(G; \frac{1}{t}, \frac{1}{1-t}\right)} & = & \sum_{A \in S(G)} k(A) t^{|A|} (1-t)^{m - |A|} - 1  \\
& = & -1 + \sum_{A \in S(G)} k(A) \sum_{j=0}^{m-|A|} (-1)^{m-|A|-j} {{m-|A|} \choose j} t^{m-j}
\eeas
This establishes the result, but we refine the coefficients further.  Define
\[p_m(t) =-1 + \sum_{A \in S(G)} k(A)
\sum_{j=0}^{m-|A|} (-1)^{m-|A|-j} {{m-|A|} \choose j} t^{m-j}.\]
Let $i=m-j$.  Then $m-|A|-j=i-|A|$, so we have
\[
p_m(t)=-1 + \sum_{A \in S(G)} k(A) \sum_{m-i=0}^{m-|A|} (-1)^{i-|A|}
{{m-|A|} \choose m-i} t^{i}.
\]

To find the coefficient of $t^i,$ we sum over all $A$ in $S(G)$ such that $|E(A)|\leq i$.  This yields 
\begin{equation}
\label{eqn:origai}
a_i= \sum_{\ell = 0}^i (-1)^{i-\ell} {{m - \ell} \choose {m-i}}  \sum_{A \in S_\ell} k(A),
\end{equation}
where $S_\ell := \{ A \in S(G) : |E(A)| = \ell \}$. Thus 
\[ p_m(t)=-1+\sum_{i=0}^ma_it^i\] with $a_i$ as above.
\ep

In the proof of Theorem \ref{thm:pmt}, we derived an initial formula (\ref{eqn:origai}) for the coefficients of the polynomial integrand in Steele's formula for the expected length of the MST.  In the next section, we derive an easier working form for the coefficients but we end this section with our first main result on the first three coefficients.

\begin{thm}
\label{thm:a0a2}
Let
\[ p_m(t)=-1+\sum_{i=0}^ma_it^i \]
be the polynomial integrand in Steele's formula for the expected length of the MST of a connected graph $G$ with $n$ vertices and $m$ edges.  Then
\[ a_0 = n, \hsp a_1 = -m, \hsp \mbox{and} \hsp a_2 = 0 \]
\end{thm}
\bp
The set $S_0$ consists of just the single subgraph of $G$ with no edges and $n$ vertices, which has $n$ connected components.  Therefore, $\sum_{A \in S_0} k(A) = n$.  Next, the set $S_1$ consists of the $m$ spanning subgraphs with just one edge, each of which has exactly $n-1$ connected components.  Therefore, $\sum_{A \in S_1} k(A) = m(n-1)$.  Lastly, the set $S_2$ consists of $m \choose 2$ spanning subgraphs with two edges, each of which has exactly $n-2$ connected components.

Substituting in these values into formula (\ref{eqn:origai}) yields
\[ a_0 = \sum_{A \in S_0} k(A) = n, \hsp a_1 = -m \sum_{A \in S_0}k(A) + \sum_{A \in S_1} k(A) = -mn + m(n-1) = -m \]
and
\beas
a_2 & =  & {m \choose 2} \sum_{A \in S_0} k(A) - (m-1)\sum_{A \in S_1}k(A) + \sum_{A \in S_2}k(A) \\
& =  & {m \choose 2} n - m (m-1) (n-1) + {m \choose 2} (n-2) = 0.
\eeas
 This completes the proof.
\ep

\subsection{Coefficients of the polynomial integrand} 

In the previous theorem, the initial formula (\ref{eqn:origai}) for the coefficients is easily applied for the cases $\ell = 0, 1,2$, because for each such $\ell$, the members of $S_\ell$ all have the same number of connected components.  When $k(A)$ is non-constant on $S_\ell,$ the enumeration becomes more difficult.

Accordingly, we partition the set $S_\ell$ into subsets with different numbers of connected components.  This can be achieved by partitioning over the ranks of the members of $S_\ell$ since subgraphs in $S_\ell$ with the same rank have the same number of connected components, namely $n-r$.

Let $k_r^\ell$ be the number of spanning subgraphs of $G$ in $S_\ell$ with rank $r$; i.e. the number of spanning subgraphs of $G$ with $\ell$ edges and $n-r$ connected components.  In terms of $k_r^\ell$, formula (\ref{eqn:origai}) can be rewritten as
\begin{equation}
\label{eqn:newai}
a_i= \sum_{\ell = 0}^i (-1)^{i-\ell} {{m - \ell} \choose {m-i}}  \sum_{r = r_\ell}^\ell k_r^\ell(n-r),
\end{equation}
where $r_\ell$ is the minimum rank of a graph with $n$ vertices and $\ell$ edges.  If $K_q$ is the largest complete graph with $|E(K_q)| < \ell$, then $r_\ell = q$.  In other words, $r_\ell$ is the largest integer with ${r_\ell \choose 2} < \ell$.

We use the fact that $\sum_{r = r_\ell}^\ell k_r^\ell = {m \choose \ell}$ to reduce the number of terms of $k_r^\ell$ by one in (\ref{eqn:newai}).  The new general expression for the polynomial coefficients $a_i$ for $i \geq 3$ is stated in Theorem \ref{thm:coeff} below.  The proof of the theorem requires a couple of combinatorial identities stated in the following lemma.

\begin{lemma}
\label{lemma:identity1} For integers $m, k, i$ and $n$,

\medskip

{\em (a)}  ${\dstyle {{m-k} \choose {m-i}} {m \choose k} = {m \choose {m-i}} {i \choose k} }$

\medskip

{\em (b)} ${\dstyle \sum_{k=0}^n(-1)^k{n\choose k}k = 0 }$
\end{lemma}

\begin{thm}
\label{thm:coeff}
Let $a_i$ be the coefficients of the polynomial integrand in Steele's integral formula for the expected length of the MST of a connected graph $G$ with $n$ vertices and $m$ edges.  Then for $i \geq 3$
\[ a_i = \sum_{\ell = 3}^i (-1)^{i-\ell} {{m - \ell} \choose {m-i}} \sum_{r = r_\ell}^{\ell -1} k_r^\ell(\ell-r).\]
\end{thm}
\bp
Summing all the terms $k_r^\ell$ for a fixed number of edges $\ell$ yields the total number of spanning subgraphs in $S_\ell$, which equals $m \choose \ell$.  This implies that $k_\ell^\ell = {m \choose \ell}  - \sum_{r = r_\ell}^{\ell - 1} k_r^\ell$ and thus from formula (\ref{eqn:newai}), we get
\beas
a_i
& = & \sum_{\ell = 0}^i (-1)^{i-\ell} {{m - \ell} \choose {m-i}} \left[ \sum_{r = r_\ell}^{\ell -1} k_r^\ell(n-r) + \left( {m \choose \ell} - \sum_{r = r_\ell}^{\ell - 1} k_r^\ell \right) (n - \ell) \right] \\
& = & \sum_{\ell = 0}^i (-1)^{i-\ell} {{m - \ell} \choose {m-i}} \left[ \sum_{r = r_\ell}^{\ell -1} k_r^\ell(\ell-r) +  {m \choose \ell} (n - \ell) \right] \\
& = & \sum_{\ell = 0}^i (-1)^{i-\ell} {{m - \ell} \choose {m-i}} \sum_{r = r_\ell}^{\ell -1} k_r^\ell(\ell-r) + \sum_{\ell = 0}^i (-1)^{i-\ell} {{m - \ell} \choose {m-i}} {m \choose \ell} (n - \ell) 
\eeas
The minimum ranks for $\ell = 0, 1, 2$ are $r_0 = 0, r_1 = 1$ and $r_2 = 2$.  Therefore, for these values of $\ell$, the summation on $r$ is empty and $a_i$ reduces to the second summation.  This and Lemma \ref{lemma:identity1}(a) yield
\[
a_i = \left[ \sum_{\ell = 3}^i (-1)^{i-\ell} {{m - \ell} \choose {m-i}} \sum_{r = r_\ell}^{\ell -1} k_r^\ell(\ell-r) \right] + \left[ \sum_{\ell = 0}^i (-1)^{i-\ell} {{m} \choose {m-i}} {i \choose \ell} (n - \ell) \right]
\]
The second sum equals zero by the Binomial Theorem and Lemma \ref{lemma:identity1}(b).
\ep

The above result gives a general formula for the coefficients of the polynomial integrand in terms of the values $k_r^\ell$.  Determining the values of $k_r^\ell$ for large $\ell$ poses a major challenge.  We conclude this section with the enumeration for $\ell = 3, 4, 5, 6$ and the 
corresponding coefficients of $p_m(t).$

\begin{defn}
For a connected graph $G$, define

{\em (a)} $c_i =$ number of cycles of size $i$ in $G$.

{\em (b)} $c_{i,1} =$ number of cycles of size $i$ with one chord.

{\em (c)} $\bar{c}_{i,1} =$ number of cycles of size $i$ with one chord and one additional edge that is not a chord of the cycle.

{\em (d)} $k_i =$ number of complete subgraphs $K_i$ in $G$.

{\em (e)} $k_{i,j} =$ number of complete bipartite subgraphs $K_{i,j}$ in $G$.

\end{defn}

\begin{lemma}
\label{lemma:dell}
For $\ell = 3, 4, 5, 6$,
\be
\label{eqn:dell}
\sum_{r = r_\ell}^{\ell - 1} k_r^\ell (\ell -r) = \sum_{j=r_\ell}^\ell c_j {{m-j} \choose {m-\ell}} - d_\ell
\ee
where
\[ d_3 = 0, \hsp d_4 = 0, \hsp d_5 = k_3^5, \hsp d_6 = \bar{c}_{4,1}+c_{5,1} + k_{3,2} + 4k_4 \]
\end{lemma}
\bp We show the above result for $\ell = 5$; the other cases are similar in nature.  The minimum rank for $\ell =5$ is $r_5 = 3$ and so the left-hand side of equation (\ref{eqn:dell}) is $k_4^5 + 2k_3^5$.  The types of subgraphs counted in $k_4^5$ are those with $5$ edges and $n-4$ connected components, which have the form shown in Figure 1(a)-(c).  Analogously, there is only one type of subgraph counted in $k_3^5$, which is shown in Figure 1(d).  Note that the graphs in Figures 1 and 2 that are a one-clique sum of smaller graphs actually
represent families that include subgraphs that are disjoint unions of the summands.  For example, 1(a) includes $K_3\bigcup P_2$, where $K_3$ is the complete graph on three vertices and $P_2$ is a path with two edges.  


\begin{figure}[h]
\begin{center}
\includegraphics[height=1.2in]{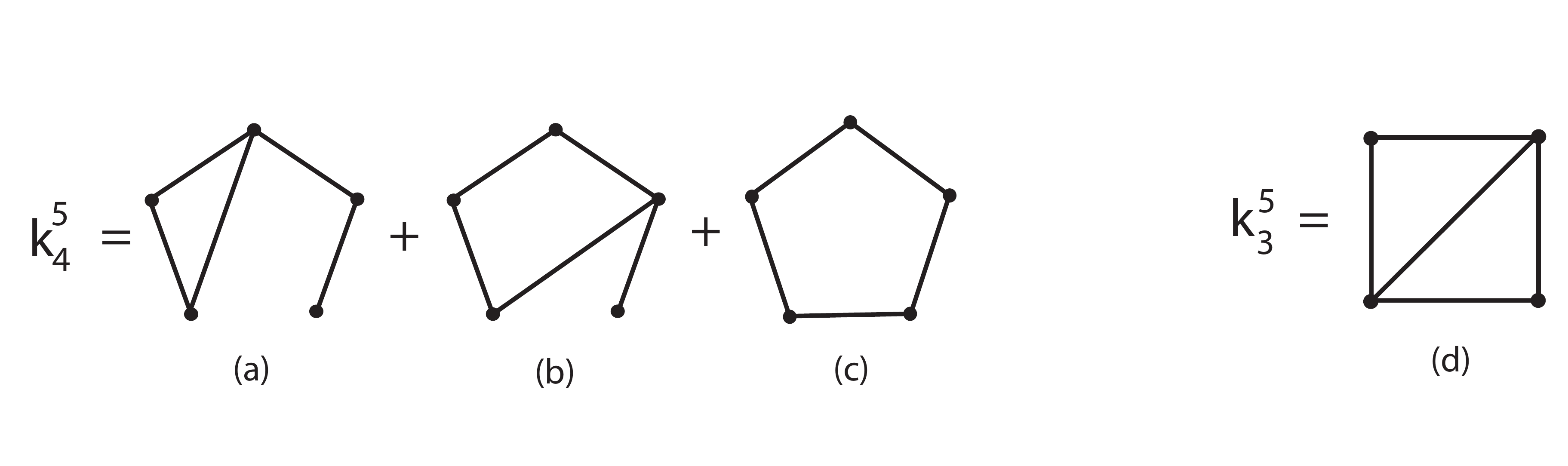}
\vspace{-.2in}
\caption{Representations of the subgraphs counted in $k_4^5$ and $k_3^5$}\label{k45andk35}
\end{center}
\end{figure}

\begin{figure}[h]
\begin{center}
\includegraphics[height=.8in]{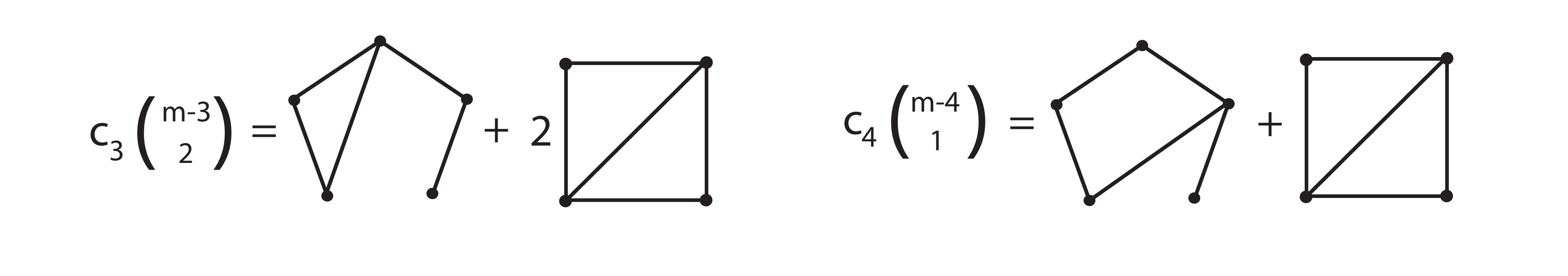}
\vspace{-.1in}
\caption{Representations of the subgraphs counted in $c_3{m-3\choose2}$ and $c_4{m-4\choose 1}$}\label{c3+andc4+}
\end{center}
\end{figure}

Now consider the right-hand side of (\ref{eqn:dell}).  Start with any 3-cycle and choose any other 2 edges in the graph; there are $c_3{{m-3} \choose 2}$ ways to do this.  This counts all the types of subgraphs depicted in Figure 1(a) and counts the subgraphs in Figure 1(d) twice.  Figure \ref{c3+andc4+} gives a pictorial representation of $c_3 {{m-3} \choose 2}.$  The subgraphs counted by $c_4{m-4\choose 1}$ (start with a 4-cycle and choose any other edge) are of the type shown in Figure \ref{k45andk35}(b) and Figure \ref{k45andk35}(d).  These are depicted in the right-hand side of Figure \ref{c3+andc4+}.

Lastly, $c_5$ is the number of $5$-cycles, which are shown in Figure \ref{k45andk35}(c).  Therefore,
\[ k_4^5 + 2k_3^5 = c_3{m-3\choose 2} + c_4{m-4\choose 1} + c_5-k_3^5. \qedhere \]
\ep

While initially Lemma \ref{lemma:dell} appears only to complicate the coefficient formula given in Theorem \ref{thm:coeff}, the next lemma shows that when it is applied to the coefficient formula, it actually simplifies it.  The proof is a straightforward calculation and so we omit it; the reasoning is analogous to the proof of Lemma \ref{lemma:dell}.  Although we proved the first equation in Lemma \ref{lemma:ci} for $i=3, 4, 5, 6$, we conjecture that it holds in general for all $i \geq 3$.

\begin{lemma}
\label{lemma:ci}
For $i = 3, 4, 5, 6$,
\[ \sum_{\ell = 3}^i (-1)^{i-\ell} {{m - \ell} \choose {m-i}} \sum_{j=r_\ell}^\ell c_j {{m-j} \choose {m-\ell}} = c_i \]
and thus
\be
\label{eqn:coeffci}
a_i = c_i - \sum_{\ell = 3}^i (-1)^{i-\ell} {{m - \ell} \choose {m-i}} d_\ell
\ee
\end{lemma}
\noi

Finally, we derive representations for the coefficients $a_3$ through $a_6$ in terms of the structure of the underlying graph $G$.  The proof of the theorem is a direct application of Lemmas \ref{lemma:dell} and \ref{lemma:ci} to the general coefficient formula given in Theorem \ref{thm:coeff}.
\begin{thm}
\label{thm:a3a6}
Let
\[ p_m(t)=-1+\sum_{i=0}^ma_it^i \]
be the polynomial integrand in Steele's formula for the expected length of the MST of a connected graph $G$ with $n$ vertices and $m$ edges.  Then
\[a_3 = c_3, \hsp a_4 = c_4, \hsp a_5 = c_5 - k_3^5, \hsp a_6 = c_6+2k_4-c_{5,1}-{k}_{3,2}. \]
\end{thm}

\section{Application of Results}\label{example} 

In this section, we apply Theorems \ref{thm:a0a2} and \ref{thm:a3a6} to the complete bipartite graph $K_{3,2}$ in order to derive the expected length of the minimal spanning tree of this graph.

\begin{proposition}
Let $p_m(t)$ be the polynomial integrand in Steele's formula for the complete bipartite graph, $K_{3,2}$ shown below.

\begin{figure}[h]
\vspace{-.20in}
\begin{center}
\includegraphics[height=.8in]{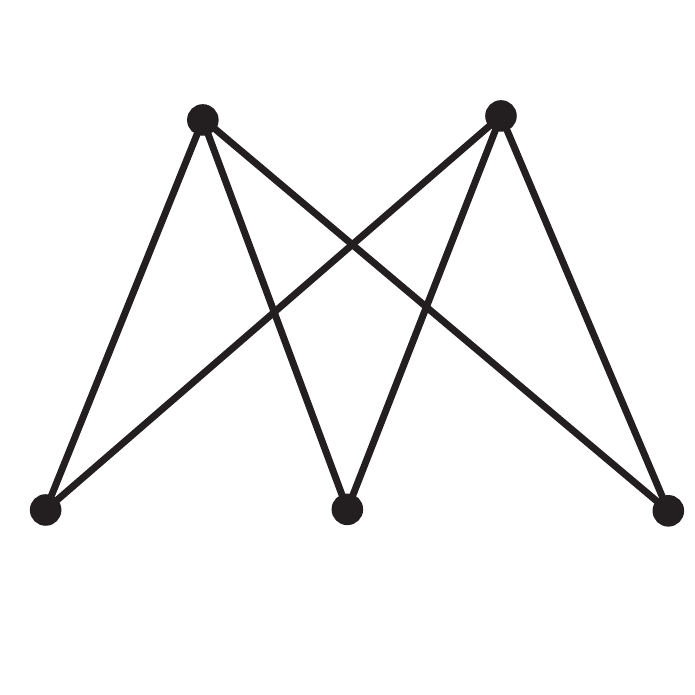}
\end{center}
\vspace{-.5in}
\end{figure}
\noi
Then
\[p_m(t) = 4 - 6t + 3t^4 -t^6\]
and 
\[ \mathbb{E}[L(K_{3,2})] = \int_0^1 (4 - 6t + 3t^4 - t^6)\, dt = 4 -3 + \frac35 - \frac17 = \frac{51}{35} \]
\end{proposition}
\bp
For $K_{3,2}$, $n = 5$, and $m=6$.  By Theorem \ref{thm:a0a2}, we have $a_0 = 5$, $a_1 = -6$, and $a_2 = 0$.  Next, we apply Theorem \ref{thm:a3a6}.  $K_{3,2}$ has no $3$-cycles, so $a_3=0$.  The graph has three 4-cycles, so $a_4 = 3$.  For the coefficient $a_5$, we note that there are no 5-cycles and also no $k_3^5$-type subgraphs (a 4-cycle with a chord) either, so $a_5 = 0$.  Lastly, for $a_6$, there are no 6-cycles, no $K_4$ subgraphs, no $c_{5,1}$-type subgraphs (since there are no 5-cycles), and one $k_{3,2}$-type subgraph (the entire graph).  Therefore, $a_6 = -1$ and we get
\[p_m(t) = -1 + 5 - 6t + 3t^4 -t^6. \qedhere \]
\ep

\section{The Complete Graph}\label{complete} 

The MST problem on $K_n$ has been studied extensively.  Frieze \cite{Frieze} proved that
\[ \lim_{n \goto \infty} \mathbb{E}[L(K_n)] = \zeta(3) = \sum_{i=1}^\infty i^{-3} = 1.202 \ldots  \]
In \cite{Ste87}, Steele extended this result to general edge distributions and Janson \cite{Janson} proved a central limit theorem for $L(K_n)$.

We apply our results to the complete graph and derive exact formulas in terms of the number of vertices $n$ for the first seven coefficients of the polynomial integrand in Steele's formula.

\begin{thm}
\label{thm:kn}
Let $p_m(t) = -1 + \sum_{i=1}^m a_i t^i$ be the polynomial integrand in Steele's formula for the complete graph on $n$ vertices, denoted by $K_n$.  Then
\[ a_0 = n, \hsp a_1 = - {n \choose 2}, \hsp a_2 = 0, \hsp a_3 = {n \choose 3} \]
\[ a_4 = 3 {n \choose 4}, \hsp a_5 = 12 {n \choose 5} - 6 {n \choose 4}, \hsp a_6 = 60 {n \choose 6} - 60 {n \choose 5} - 2 (n-5){n \choose 4} \]
\end{thm}
\bp
For the complete graph on $n$ vertices, the number of edges $m = {n \choose 2}$ and the number of cycles of length $j$ is given by
\[ c_j = \frac{1}{2} {n \choose j} (j-1)! \]
In addition, $k_3^5 = 2 c_4$, $c_{5,1} = 5c_5$, $k_4 = {n \choose 4}$ and $k_{3,2} = {n \choose 5}{5 \choose 2}$.
\ep

Numerical calculation of $\mathbb{E}[L(K_n)]$ had led to the famous conjecture that the convergent sequence is also monotone increasing and concave.  This problem was raised at the conference {\it Mathematics and Computer Science II} at Versailles in 2002 but no proof has been found.  Clearly, our results alone will not answer this question as we have only derived exact formulas for the first seven coefficients.  But our results give a hint that there may be a pattern to the coefficients of the polynomial integrand in Steele's formula for the complete graph, which if true, would answer the conjecture.

We end this section with a result that factors the polynomial integrand in Steele's formula for $K_n$, with one of the factors a polynomial of degree less than or equal to the number of edges of $K_{n-1}$.

\begin{thm}
Let $p_m(t)$ be the polynomial integrand in Steele's formula for the expected length of the MST of the complete graph on $n$ vertices denoted by $K_n$.  Then
\[ p_m(t) = (1-t)^{n-1} q (t), \]
where $q(t)$ is a polynomial of degree less than or equal to ${{n-1} \choose 2}$.
\end{thm}
\bp
As in the proof of Theorem \ref{thm:pmt}, we have \[p_m(t)=\sum_{A\in S(G)}k(A)t^{|A|}(1-t)^{m-|A|}-1,\] where $m={n\choose 2}$.  

We factor out $(1-t)^{n-1}$ to get
\[
p_m(t)= (1-t)^{n-1}\left[\sum_{A\in S(G)}k(A)t^{|A|}(1-t)^{{n\choose 2}-|A| - (n-1)}-(1-t)^{1-n}\right].\]  Note that ${n\choose 2} - (n-1) = {n-1\choose 2}$.
Now the sum ranges over spanning subgraphs of size (in edges) from 0 to ${n\choose 2}$.  We split it into two sums as follows:
\beas
\lefteqn{ p_m(t)=(1-t)^{n-1}\left[\sum_{|A|\leq {n-1\choose 2}} k(A)t^{|A|}(1-t)^{{n - 1\choose 2} - |A|} \right.} \\
& + & \left. \sum_{|A|>{n-1\choose 2}}k(A)t^{|A|}(1-t)^{{n \choose 2} - |A| - (n-1)} - (1-t)^{1-n}\right]
\eeas
Clearly, the first sum over $|A|\leq {n-1\choose 2}$ is a polynomial of degree at most ${n-1\choose 2}$.  Call it $q_1(t)$.

For the second sum, we sum over possible number of edges $i> {n-1\choose 2}$ and count the number of subgraphs with $i$ edges, which for the complete graph is ${{n\choose 2} \choose i}$.  Furthermore,  for all spanning subgraphs of $K_n$ with $i> {n-1\choose 2}$ edges, the number of connected components is $1$. Therefore, we have

\beas
\lefteqn{\sum_{|A|>{n-1\choose 2}}k(A)t^{|A|}(1-t)^{{n \choose 2} - |A| - (n-1)} } \\
& = & \sum_{i= {n-1\choose 2} + 1}^{n\choose 2}{{n\choose 2} \choose i}t^i(1-t)^{{n\choose 2} -i - (n-1)} \\
& = & \sum_{i=0}^{n\choose2} {{n\choose 2}\choose i}t^i(1-t)^{{n\choose 2}-i}(1-t)^{-(n-1)} - \sum_{i=0}^{{n-1\choose2}}{{n\choose 2}\choose i}t^i(1-t)^{{n - 1\choose 2}-i} \\
& = & (1-t)^{1-n}\sum_{i=0}^{n\choose2} {{n\choose 2}\choose i}t^i(1-t)^{{n\choose 2}-i} - \sum_{i=0}^{{n-1\choose2}}{{n\choose 2}\choose i}t^i(1-t)^{{n - 1\choose 2}-i}
\eeas
By the Binomial Theorem, the first sum equals $1$ and the second sum, call it $q_2(t),$ is a polynomial of degree at most ${n-1\choose 2}$.

We now have
\[p_m(t)=(1-t)^{n-1}(q_1(t) + (1-t)^{1-n} + q_2(t) - (1-t)^{1-n})=(1-t)^{n-1}(q_1(t) + q_2(t)),\] where both $q_1(t)$ and $q_2(t)$ are polynomials of degree less than or equal to ${n-1\choose 2}$.  This completes the proof.
\ep

\section*{About the authors:}

\subsection*{Jared Nishikawa}
   Jared Nishikawa graduated with a B.A. in Mathematics from Willamette University in May '10.
   He is currently in graduate school at the University of Colorado at Boulder, pursuing a Ph.D. in Mathematics.

   1819 NE 205th Ave., Fairview, OR 97024, jared.nishikawa@gmail.com

\subsection*{Peter T. Otto}

Peter T. Otto earned his Ph.D. in mathematics from the University of Massachusetts, Amherst in 2004.  His primary research interest is in probability theory with a focus on the study of models in statistical mechanics.

    Willamette University, Salem, Oregon 97301,  potto@willamette.edu

\subsection*{Colin Starr}

Colin Starr earned his Ph.D. in mathematics from the University of Kentucky in 1998.  His primary research interest is in graph theory.

 Willamette University, Salem, Oregon 97301, cstarr@willamette.edu

\end{document}